\documentclass[a4paper,12pt]{amsart}
\usepackage{enumitem}
\usepackage{bigints}
\usepackage{amssymb, amsmath, mathtools}
\usepackage{graphicx}
\usepackage{ifthen}
\usepackage{xcolor}
\usepackage{cite}
\nonstopmode \numberwithin{equation}{section}
\newtheorem{thm}{Theorem}
\newtheorem{cor}{Corollary}
\newtheorem{lem}{Lemma}
\newtheorem{prop}{Proposition}

\newtheorem{conj}{Conjecture}

\theoremstyle{definition}
\newtheorem{defn}{Definition}[section]

\newtheorem{prob}[equation]{Problem}

\newenvironment{rem}{%
\bigskip
\noindent \textsl{{\bf Remark. }}}{\bigskip}
\newenvironment{rems}{%
\bigskip
\noindent \textsl{{\bf Remarks. }}}{\bigskip}

\newcounter {own}
\def\theown {\thesection       .\arabic{own}}

\newenvironment{pf}[1][]{%
 \vskip 3mm
 \noindent
 \ifthenelse{\equal{#1}{}}%
  {{\slshape Proof. }}%
  {{\slshape #1.} }%
 }%
{\qed\bigskip}

\newcounter{alphabet}

\setlist[itemize]{leftmargin=*}

\newcommand{\ID}{{\mathbb D}}
\newcommand{\IN}{{\mathbb N}}
\newcommand{\IC}{{\mathbb C}}

\def\be{\begin{equation}}
\def\ee{\end{equation}}

\newcommand{\bee}{\begin{enumerate}}
\newcommand{\eee}{\end{enumerate}}

\newcommand{\blem}{\begin{lem}}
\newcommand{\elem}{\end{lem}}
\newcommand{\bthm}{\begin{thm}}
\newcommand{\ethm}{\end{thm}}
\newcommand{\bcor}{\begin{cor}}
\newcommand{\ecor}{\end{cor}}
\newcommand{\beg}{\begin{examp}}
\newcommand{\eeg}{\end{examp}}
\newcommand{\begs}{\begin{examples}}
\newcommand{\eegs}{\end{examples}}
\newcommand{\bdefe}{\begin{defn}}
\newcommand{\edefe}{\end{defn}}
\newcommand{\bprob}{\begin{prob}}
\newcommand{\eprob}{\end{prob}}
\newcommand{\bei}{\begin{itemize}}
\newcommand{\eei}{\end{itemize}}

\newcommand{\bcon}{\begin{conj}}
\newcommand{\econ}{\end{conj}}
\newcommand{\bcons}{\begin{conjs}}
\newcommand{\econs}{\end{conjs}}
\newcommand{\bprop}{\begin{prop}}
\newcommand{\eprop}{\end{prop}}
\newcommand{\br}{\begin{rem}}
\newcommand{\er}{\end{rem}}
\newcommand{\brs}{\begin{rems}}
\newcommand{\ers}{\end{rems}}
\newcommand{\bo}{\begin{obser}}
\newcommand{\eo}{\end{obser}}
\newcommand{\bos}{\begin{obsers}}
\newcommand{\eos}{\end{obsers}}
\newcommand{\bpf}{\begin{pf}}
\newcommand{\epf}{\end{pf}}
\newcommand{\ba}{\begin{array}}
\newcommand{\ea}{\end{array}}
\newcommand{\beq}{\begin{eqnarray}}
\newcommand{\beqq}{\begin{eqnarray*}}
\newcommand{\eeq}{\end{eqnarray}}
\newcommand{\eeqq}{\end{eqnarray*}}

\newcommand{\ov}{\overline}

\newcounter{minutes}
\divide\time by 60
\newcounter{hours}
\multiply\time by 60 \addtocounter{minutes}{-\time}
\begin{document}
\title{The $p$-Bohr radius for vector-valued holomorphic and pluriharmonic functions}
\author{Nilanjan Das}
\address{Nilanjan Das, Theoretical Statistics and Mathematics Unit, Indian Statistical Institute Kolkata, Kolkata-700108, India.}
\email{nilanjand7@gmail.com}


\begin{abstract}
We study a ``$p$-powered" version $K_n^p(F(R))$ of the well-known Bohr radius problem for the family $F(R)$ of holomorphic functions $f: R\to X$ satisfying $\|f\|<\infty$, where $\|.\|$ is a norm in the function space $F(R)$, $R\subset\IC^n$ is a complete Reinhardt domain and $X$ is a complex Banach space. For all $p>0$, we describe in full details the asymptotic behaviour of $K_n^p(F(R))$, where 
$F(R)$ is (a) the Hardy space of $X$-valued holomorphic functions defined in the open unit polydisk $\ID^n$, and (b) the space of bounded $X$-valued holomorphic or complex-valued pluriharmonic functions defined in the open unit ball $B(l_t^n)$ of the Minkowski space $l_t^n$. We give an alternative definition of the optimal cotype for a complex Banach space $X$ in the light of these results.
In addition, the best possible versions of two theorems from [B\'en\'eteau et. al., Comput. Methods Funct. Theory, 4 (2004), no. 1, 1--19] and [Chen \& Hamada, J. Funct. Anal., 282 (2022), no. 1, Paper No. 109254, 42 pp] have been obtained as specific instances of our results.
\end{abstract}
\thanks{The author of this article is supported by a Research Associateship provided by the Stat-Math Unit of ISI Kolkata.}

\maketitle
\pagestyle{myheadings}
\markboth{The $p$-Bohr radius for vector-valued holomorphic and pluriharmonic functions}{N. Das}

\bigskip
\section{Introduction and the main results}\label{N3sec1}
A complete Reinhardt domain $R$ in the $n$-dimensional complex plane $\IC^n$ is a domain such that if $z=(z_1, z_2,\cdots, z_n)\in R$ then $x=(x_1, x_2, \cdots, x_n)\in R$, provided that $|x_n|\leq |z_n|$ for all $n\in\IN$. Throughout this article, we will use the standard multi-index notation: $\alpha$ denotes an $n$-tuple $(\alpha_1, \alpha_2,\cdots, \alpha_n)$ of
nonnegative integers, $|\alpha|:=\alpha_1+\alpha_2+\cdots+\alpha_n$, $\alpha!:=\alpha_1!\alpha_2!\cdots\alpha_n!$, $z$ denotes an $n$-tuple $(z_1, z_2, \cdots, z_n)$ of
complex numbers, and $z^\alpha$ is the product $z_1^{\alpha_1}z_2^{\alpha_2}\cdots z_n^{\alpha_n}$. 
Given a complex Banach space $X$, let us define the \emph{Bohr radius} $K_X(R)$ as the supremum of all $r\in[0,1]$ such that
\be\label{N3eq19}
\sup_{z\in rR} \sum_{k=0}^\infty\sum_{|\alpha|=k}\|x_\alpha z^\alpha\|\leq 1
\ee
for all holomorphic functions $f(z)=\sum_{k=0}^\infty\sum_{|\alpha|=k}x_\alpha z^\alpha:R\to X$ with $\|f(z)\|\leq 1$ for all $z\in R$.
A classical theorem of Harald Bohr (optimized independently by Wiener, Riesz and Schur) states that $K_\IC(\ID)=1/3$ (cf. \cite{Bohr}). In 1995, Dixon \cite{Dix} found an unexpected application of this result in the realm of operator algebras. Since then, Bohr's theorem has been studied and extended from several different aspects-see, for example, the articles \cite{Aiz, Bala, Bay, Def2, Seip, Ga, Ham, La, Li, Paul, Pop1}, the recent book \cite{Def3} and their references.

In this article, we will concentrate on the $n$-dimensional generalizations and the vector-valued aspects of this theory. The quantity $K_\IC(\ID^n)$ was considered for the first time by Dineen and Timoney \cite{Din}, and subsequently by Boas and Khavinson \cite{Boas}, where $\ID^n:=\{z=(z_1, z_2, \cdots, z_n)\in\IC^n:|z_i|<1\,\textrm{for all}\, 1\leq i\leq n\}$ is the open unit polydisk in $\IC^n$. It is worth mentioning here that the exact value for $K_\IC(\ID^n)$ for all $n>1$ remains unknown till now. This makes the study of the asymptotic behaviour of $K_\IC(\ID^n)$ as a function of $n$ an intriguing problem. The best known result in this direction is the following one by Bayart et. al. \cite{Bay}:
$$
\lim_{n\to\infty}\frac{K_\IC(\ID^n)}{\sqrt{\frac{\log n}{n}}}=1
$$
(see also \cite{Seip}). Another important contribution is the undermentioned theorem of Defant and Frerick \cite{Def6}:
$$
K_\IC(B(l_t^n))\sim\left(\frac{\log n}{n}\right)^{1-\frac{1}{\min\{t, 2\}}},
$$
where $B(l_t^n):=\{z=(z_1, z_2, \cdots, z_n)\in\IC^n:\left(\sum_{i=1}^n|z_i|^t\right)^{1/t}<1\}$ is the open unit ball of the Minkowski space $l_t^n$, $t\in[1, \infty)$, and $\ID^n\equiv B(l_\infty^n)$. We clarify here that for any two sequences $\{p_n\}$ and $\{q_n\}$ of positive real numbers, we write $p_n\prec q_n$ and $p_n\succ q_n$ if there exist constants $C, D>0$ such that
$p_n\leq Cq_n$ and $p_n\geq Dq_n$ for all $n>1$, respectively; and $p_n\sim q_n$ whenever $p_n\prec q_n$ and $p_n\succ q_n$. 

At this point, it seems appropriate to mention that the classical form $(\ref{N3eq19})$ of the Bohr inequality becomes unusable for $\dim(X)>1$, as we already know $K_{l_t^2}(\ID)=0$ for all $t\in[1, \infty]$ (see \cite[Theorem 1.2]{Bla1}). Yet, proper modification of this inequality generates interesting new results involving the corresponding ``Bohr radius". It was shown in \cite{Def} by Defant et. al. that if we replace $1$ in the right hand side of the inequality (\ref{N3eq19}) by $\lambda>1$
and define in a similar way the Bohr radius $K_X(\ID^n, \lambda)$ for this case, then $K_X(\ID^n, \lambda)>0$ for all complex Banach spaces $X$. Moreover, \cite[Theorem 1.2]{Def} includes the precise asymptotic behaviour of $K_X(\ID^n, \lambda)$ for both finite and infinite dimensional Banach spaces.

On the other hand, if one chooses to keep the right hand side of $(\ref{N3eq19})$ intact and reshape the left hand side, then a way to do so is to take sum of the positive powers of the terms $\|x_\alpha z^\alpha\|$. Namely, for $p>0$ let us consider the Bohr inequality
$$\sup_{z\in r\ID^n} \sum_{k=0}^\infty\sum_{|\alpha|=k}\|x_\alpha z^\alpha\|^p\leq 1$$
for all holomorphic functions $f(z)=\sum_{k=0}^\infty\sum_{|\alpha|=k}x_\alpha z^\alpha:\ID^n\to X$ with $\|f(z)\|\leq 1$ in $\ID^n$, and denote by $K_X^p(\ID^n)$ the largest possible $r\in [0,1]$ for which the above inequality remains valid. 

The quantities $K_\IC^p(\ID)$ and $K_\IC^p(\ID^n)$ were first studied by Djakov and Ramanujan \cite{Dja}. It is known from their work that both of these quantities are $0$ for $p\in(0,1)$ and are equal to $1$ for $p\geq 2$. Analogues of these radii in Hardy spaces were investigated by B\'en\'eteau et. al. \cite{Bene}. For $p\in(1,2)$, the exact value of $K_\IC^p(\ID)$ was determined by Kayumov and Ponnusamy \cite{Kay1}, and the exact asymptotic value of $K_\IC^p(\ID^n)$ has recently been obtained by the author \cite{Das}.
For vector-valued situations, Blasco showed in \cite[Theorem 1.10]{Bla2} that $K_X^p(\ID)>0$ if and only if $X$ is \emph{$p$-uniformly $\IC$-convex} ($2\leq p<\infty$), and that $K_X^p(\ID)$ cannot have nonzero value when $p\in[1, 2)$ and $\dim(X)\geq 2$. Thus, such a definition does not produce a nonzero Bohr radius for all complex Banach spaces $X$. This raises the following 

\vspace{3pt}
\noindent
\textbf{\underline{Problem:}} \emph{To get a version of the above-mentioned powered Bohr inequality that yields a nonzero Bohr radius for any complex Banach space $X$.}

\vspace{3pt}
\noindent
We propose a slightly different definition to tackle this question.
Let us denote by $F(R)$ the family of holomorphic functions $f: R\to X$ satisfying $\|f\|<\infty$, $\|.\|$ being a norm in the function space $F(R)$. Here $R\subset\IC^n$ is a complete Reinhardt domain and $X$ is a complex Banach space. For any given $p\in(0, \infty)$, the $p$-Bohr radius $K_n^p(F(R))$ for the family $F(R)$ is the supremum of all $r\in [0, 1]$ such that 
$$
\sup_{z\in rR} \sum_{k=1}^\infty\sum_{|\alpha|=k}\|x_\alpha z^\alpha\|^p\leq 1
$$
for all $f(z)=\sum_{k=0}^\infty\sum_{|\alpha|=k}x_\alpha z^\alpha\in F(R)$ satisfying $\|f\|\leq 1$. 

It will soon appear as a special case of our results that for all $p>0$, $K_n^p(F(R))>0$ for any complex Banach space $X$ and any complete Reinhardt domain $R\subset\IC^n$, provided that $F(R)$ is the family of bounded (w.r.t the sup-norm) $X$-valued holomorphic functions defined in $R$. 

We will now start presenting our main theorems, along with necessary background information. To the best of our knowledge, after the article \cite{Bene} there has not been any progress on the study of the powered Bohr inequality for Hardy space functions in several variables. This motivates us to consider the following 

\vspace{3pt}
\noindent
\textbf{\underline{Problem:}} \emph{To study in detail the quantity $K_n^p(F(R))$, where $F(R)$ is the Hardy space of $X$-valued holomorphic functions defined in $\ID^n$.}

\vspace{3pt}
\noindent
We will address this question in Theorem \ref{N3thm1}, but at first we need to introduce the notion of vector-valued Hardy space and the ``cotype" of a Banach space $X$. Namely, for $s\in[1, \infty]$ we define
$$
H^s(\ID^n, X):=\{f:\ID^n\to X\,\,\mbox{holomorphic}\,:\|f\|_s<\infty\},
$$
where
$$
\|f\|_s:=
\begin{cases}
\sup_{0<\rho<1}\left(\int_{\mathbb{T}^n}\|f(\rho z)\|^s dz\right)^{\frac{1}{s}}\,\,\mbox{for}\,\, s<\infty,\,\,\mbox{and}\\
\sup_{z\in\ID^n}\|f(z)\|\,\,\mbox{for}\,\, s=\infty,
\end{cases}
$$
$dz$ being the normalized Lebesgue measure on the unit polytorus $\mathbb{T}^n:=\{z=(z_1, z_2, \cdots, z_n)\in\IC^n:|z_i|=1\,\textrm{for all}\, 1\leq i\leq n\}$. Furthermore, a Banach space $X$ is said to have cotype $q\in[2, \infty]$ if there exists a constant $C>0$ such that for arbitrarily chosen vectors $x_1, x_2,\cdots, x_N\in X$ we have
$$
\left(\sum_{n=1}^N\|x_n\|^q\right)^{\frac{1}{q}}\leq C\left(\int_0^1\left\|\sum_{n=1}^Nr_n(t)x_n\right\|^2dt\right)^{\frac{1}{2}},
$$
where $r_n$ is the $n$-th Rademacher function on $[0, 1]$. We write
$$
\textrm{Cot}(X):=\inf\{2\leq q\leq\infty: X\,\,\mbox{has cotype}\,\, q\}.
$$
The quantity $\mbox{Cot}(X)$ is sometimes referred to as the \lq\lq optimal cotype\rq\rq of $X$, and the infimum in its definition may not be attained. The concept of cotype plays an important role in the study of summability theory and geometry of Banach spaces. We urge the interested reader to look into the references \cite{Die, Def3} for more information on these topics.

With all these preparations, we are now ready to state the first theorem of this article.
\bthm\label{N3thm1}
For any given complex Banach space $X$ and $p\in (0, \infty), s\in[1, \infty]$, we have the following:
\begin{itemize}
\item[\emph{(i)}]
When $\dim(X)=\infty$, 
\begin{align*}\left(\frac{1}{n}\right)^{\frac{1}{p}-\frac{1}{q}}&\prec K_n^p(H^s(\ID^n, X))\prec \left(\frac{1}{n}\right)^{\frac{1}{p}-\frac{1}{\emph{Cot}(X)}}\,\emph{for}\,p\leq\emph{Cot}(X),\,\emph{and}\\
& K_n^p(H^s(\ID^n, X))\sim 1\,\emph{for}\, p>\emph{Cot}(X),
\end{align*}
if $X$ has a finite cotype $q$. Otherwise
$$
K_n^p(H^s(\ID^n, X))\sim\left(\frac{1}{n}\right)^{\frac{1}{p}}.
$$
\item[\emph{(ii)}]
When $\dim(X)<\infty$ and $p<2$,
$$
K_n^p(H^s(\ID^n, X))\sim
\begin{cases}
\left(\frac{1}{n}\right)^{\frac{1}{p}-\frac{1}{2}}\,\emph{for}\,s\in[1,2],\,\emph{and}\\
\left(\frac{\log n}{n}\right)^{\frac{1}{p}-\frac{1}{2}}\,\emph{for}\, s\in(2, \infty],
\end{cases}
$$
and $K_n^p(H^s(\ID^n, X))\sim 1$ for all $s\geq 1, p\geq 2$.
\item[\emph{(iii)}]
In particular, for all $p<2$
$$
K_n^p(H^2(\ID^n, \IC))=\left(1-\left(\frac{1}{2}\right)^{\frac{1}{n}}\right)^{\frac{1}{p}-\frac{1}{2}}.
$$
\end{itemize}
\ethm
\br
It is interesting to note that for $p=1$ and $s=\infty$, part (i) and part (ii) of this theorem coincide with \cite[Theorem 1.2]{Def}. Also, part (ii) coincides with \cite[Theorem 2]{Seip} for $X=\IC$ and $p=1,\,s=\infty$.
\er

The following corollary is immediate from part (iii) of the Theorem \ref{N3thm1} (see the proof also).
\bcor\label{N3cor1}
For any $p\in(0,2)$ and for any $f(z)=\sum_{k=1}^\infty\sum_{|\alpha|=k}c_\alpha z^\alpha\in H^2(\ID^n, \IC)$ with $\|f\|_2\leq 1$, the inequality $\sum_{k=1}^\infty\sum_{|\alpha|=k}|c_\alpha|^p |z^\alpha|\leq 1$ holds for all $z\in R_n\ID^n$, where $R_n\leq (1-(1/2)^{1/n})^{1-(p/2)}$. This upper bound for $R_n$ is the best possible.
\ecor
This Corollary \ref{N3cor1} sharpens \cite[Theorem 4.1]{Bene}, where the above result was shown to hold for $R_n\leq (1/(2n))^{1-(p/2)}$. Another quick outcome of Theorem \ref{N3thm1} and Proposition \ref{N3prop1} (see section \ref{N3sec2}) is the following:
\bcor\label{N3cor2}
Given any complete Reinhardt domain $R\subset\IC^n$, $K_n^p(H^\infty(R, X))>0$ for all complex Banach spaces $X$ and for all $p\in (0, \infty)$. Here $H^\infty(R, X)$ is the space of $X$-valued holomorphic functions defined in $R$, bounded w.r.t the sup-norm.
\ecor
We now turn our attention to the spaces $H^\infty(B(l_t^n), X)$, $t\in[1, \infty]$. As we have already mentioned, exact asymptotic value for $K_n^1(H^\infty(B(l_t^n), \IC))$ is known from \cite{Def6}, while a similar quantity has been studied in \cite{Def} for functions in $H^\infty(\ID^n, X)$. The goal of the following theorem is to provide a unification of these two approaches in a $p$-powered framework.

\bthm\label{N3thm2}
For any given complex Banach space $X$ and $p\in (0, \infty), t\in[1, \infty]$, we have the following:
\begin{itemize}
\item[\emph{(i)}]
When $\dim(X)=\infty$, 
\begin{align*}\left(\frac{1}{n}\right)^{\frac{1}{p}-\frac{1}{\min\{q, t\}}}&\prec K_n^p(H^\infty(B(l_t^n), X))\prec \left(\frac{1}{n}\right)^{\frac{1}{p}-\frac{1}{\min\{\emph{Cot}(X), t\}}}\,\emph{for}\,p\leq\min\{\emph{Cot}(X), t\}\,\emph{and}\\
 & K_n^p(H^\infty(B(l_t^n), X))\sim 1\,\emph{for}\, p>\min\{\emph{Cot}(X), t\},
\end{align*}
if $X$ has a finite cotype $q$. Otherwise
$$
K_n^p(H^\infty(B(l_t^n), X))\sim\left(\frac{1}{n}\right)^{\frac{1}{p}-\frac{1}{t}}.
$$
\item[\emph{(ii)}]
When $\dim(X)<\infty$,
$$
K_n^p(H^\infty(B(l_t^n), X))\sim
\begin{cases}
\left(\frac{\log n}{n}\right)^{\frac{1}{p}-\frac{1}{\min\{t,2\}}}\,\emph{for}\, p<\min\{t,2\},\,\emph{and}\\
1\,\,\emph{for}\,p\geq\min\{t,2\}.
\end{cases}
$$
\end{itemize}
\ethm
Let us note from Theorem \ref{N3thm1} that for any chosen $s\in[1, \infty]$, if $p<\textrm{Cot}(X)$ then $K_n^p(H^s(\ID^n, X))$ tends to $0$ as $n\to\infty$. On the other hand, $K_n^p(H^s(\ID^n, X))$ remains bounded away from $0$ (as a function of $n$) for $p>\textrm{Cot}(X)$, and also for $p=\textrm{Cot}(X)<\infty$ if $X$ has $\textrm{Cot}(X)$ as cotype. Similarly, Theorem \ref{N3thm2} implies that whenever $p<\min\{\textrm{Cot}(X), t\}$, $K_n^p(H^\infty(B(l_t^n), X))$ approaches $0$ as $n\to\infty$ and for $p>\min\{\textrm{Cot}(X), t\}$, $K_n^p(H^\infty(B(l_t^n), X))$ is bounded away from $0$. Summarizing these observations, we pose the following characterization of optimal cotype-i.e. $\textrm{Cot}(X)$-of $X$.
\bcor\label{N3cor4}
For any complex Banach space $X$, we have
\begin{itemize}
\item[\emph{(i)}]
$\emph{Cot}(X)=\sup\{p>0:\lim_{n\to\infty}K_n^p(H^s(\ID^n, X))=0\}$ for any $s\in[1, \infty]$ of our choice.
\item[\emph{(ii)}]
$X$ cannot have $\emph{Cot}(X)$ as a finite cotype, provided that $K_n^{\emph{Cot}(X)}(H^s(\ID^n, X))\to 0$ as $n\to\infty$.
\item[\emph{(iii)}]
The supremum of all $p>0$ such that $\lim_{n\to\infty}K_n^p(H^\infty(B(l_t^n), X))=0\,\,(t\in[1, \infty])$ is either $t$ or $\emph{Cot}(X)$.
\end{itemize}
\ecor

Finally, we focus on complex-valued pluriharmonic functions defined in $B(l_t^n)$. A twice continuously differentiable function $f:\Omega\subset\IC^n\to\IC$ is said to be \emph{pluriharmonic} if its restriction to every complex line is harmonic. Moreover, if $\Omega$ is a simply connected domain in $\IC^n$ containing the origin, then $f:\Omega\to\IC$ is pluriharmonic if and only if it admits a representation $f=h+\ov{g}$, where $h, g$ are holomorphic in $\Omega$ with $g(0)=0$ (see \cite{Vla}). Let us now denote by $R_p(B(l_t^n))$ the supremum of all $r\in[0, 1]$ such that the inequality
$$
\sup_{z\in rB(l_t^n)} \sum_{k=1}^\infty\sum_{|\alpha|=k}\left(|a_\alpha |^p+|b_\alpha|^p\right)|z^\alpha|^p\leq 1
$$
holds for all pluriharmonic functions $f(z)=\sum_{k=0}^\infty\sum_{|\alpha|=k}a_\alpha z^\alpha+\ov{\sum_{k=1}^\infty\sum_{|\alpha|=k}b_\alpha z^\alpha}$ from $B(l_t^n)$ to $\ov{\ID}$. The following theorem shows that $R_p(B(l_t^n))$ exhibits the same asymptotic behaviour as $K_n^p(H^\infty(B(l_t^n), \IC))$.
\bthm\label{N3thm3}
For any given $p\in(0, \infty)$ and $t\in[1, \infty]$, we have 
$$
R_p(B(l_t^n))\sim
\begin{cases}
\left(\frac{\log n}{n}\right)^{\frac{1}{p}-\frac{1}{\min\{t,2\}}}\,\emph{for}\, p<\min\{t,2\},\,\emph{and}\\
1\,\,\emph{for}\,p\geq\min\{t,2\}.    
\end{cases}
$$
\ethm
The quantity $R_1(B(l_t^n))$ was introduced for the first time by Chen and Hamada \cite{Chen}. In particular, \cite[Theorem 2.11]{Chen} states that
$$
\left(\frac{1}{n}\right)^{1-\frac{1}{\min\{t, 2\}}}\prec R_1(B(l_t^n))\prec\left(\frac{\log n}{n}\right)^{1-\frac{1}{\min\{t, 2\}}}
$$
for all $t\in[1, \infty]$. For $p=1$, our Theorem \ref{N3thm3} gives
the following optimal version of this result.
\bcor\label{N3cor3}
For all $t\in[1, \infty]$, we have
$$
R_1(B(l_t^n))\sim\left(\frac{\log n}{n}\right)^{1-\frac{1}{\min\{t, 2\}}}.
$$
\ecor
\section{Preliminaries}\label{N3sec2}
This section is intended for listing some known, or unknown but easily deducible facts, which are essential for proving the three theorems stated in previous section.
Let us denote by $L^\beta(W, d\omega)\,(1\leq \beta\leq\infty)$ the Lebesgue space of complex-valued $\omega$-measurable functions $f$ defined on $W$ such that 
$$
\|f\|_{L^\beta}=
\begin{cases}
\left(\int_W|f(x)|^\beta d\omega\right)^{\frac{1}{\beta}}\,\,\mbox{for}\,\,\beta<\infty,\\
\sup_{x\in W}|f(x)|\,\,\mbox{for}\,\,\beta=\infty
\end{cases}
$$
is finite, where $\omega$ is always a positive measure. Also, we denote by $T$ a linear mapping between these spaces. With these notations, we now state the well-known Riesz-Thorin interpolation theorem as follows:

\vspace{3pt}
\noindent
\textbf{Theorem A (see f.i. \cite[Theorem 1.1.1]{Ber}).}
Suppose $p_0\neq p_1$ and $q_0\neq q_1$, and $T:L^{p_0}(U, d\mu)\to L^{q_0}(V, d\nu)$ with norm $M_0$, $T:L^{p_1}(U, d\mu)\to L^{q_1}(V, d\nu)$ with norm $M_1$. Then $T:L^{p_2}(U, d\mu)\to L^{q_2}(V, d\nu)$ is a linear mapping with norm $M_2$ satisfying
$$
M_2\leq M_0^{1-\theta}M_1^\theta,
$$
provided that $0<\theta<1$, and that
$$
\frac{1}{p_2}=\frac{1-\theta}{p_0}+\frac{\theta}{p_1},\,\, 
\frac{1}{q_2}=\frac{1-\theta}{q_0}+\frac{\theta}{q_1}.
$$

\vspace{3pt}
For $1\leq t<\infty$ and for a linear operator $A:X_0\to Y_0$ between the complex Banach spaces $X_0$ and $Y_0$, we say that $A$ is \emph{$t$-summing} if there exists a constant $c\geq 0$ such that regardless of the natural number $m$ and regardless of the choice of $f_1, f_2,\cdots, f_m$ in $X_0$, we have
$$
\left(\sum_{i=1}^m\|A(f_i)\|^t\right)^{1/t}\leq c\sup_{\phi\in B(X_0^*)}\left(\sum_{i=1}^m|\phi(f_i)|^t\right)^{1/t},
$$
where $B(X_0^*)$ is the open unit ball in the dual space $X_0^*$. The least $c$ for which the above inequality always holds is denoted by $\pi_t(A)$, and the set of all $t$-summing operators from $X_0$ into $Y_0$ is denoted by
$\Pi_t(X_0, Y_0)$. Now, the following facts are known:

\vspace{3pt}
\noindent
\textbf{Theorem B (see \cite[Proposition 2.3]{Die}).}
If $A:X_0\to Y_0$ is a bounded linear operator and $\mbox{dim}(A(X_0))<\infty$, then $A$ is $t$-summing for every $t\in[1, \infty)$.

\vspace{3pt}
\noindent
\textbf{Theorem C (see \cite[Theorem 2.8]{Die}).}
If $1\leq t_1<t_2<\infty$, then
$\Pi_{t_1}(X_0, Y_0)\subset\Pi_{t_2}(X_0, Y_0)$. Moreover, for $A\in\Pi_{t_1}(X_0, Y_0)$, we have $\pi_{t_2}(A)\leq \pi_{t_1}(A)$.

\vspace{3pt}
We also recall the following result from \cite{Bay1}, which is instrumental in obtaining the upper bound on $K_n^p(H^\infty(B(l_t^n), \IC))$ (and therefore on $K_n^p(H^\infty(B(l_t^n), X))$ with $\dim(X)<\infty$ as well).

\vspace{3pt}
\noindent
\textbf{Theorem D (cf. \cite[Corollary 3.2]{Bay1}).} Let $Y=(\IC^n, \|.\|_Y)$ be an $n$-dimensional complex Banach space with $B(Y)$ its open unit ball and $k\geq 2$. Let also $u\in[1, 2]$ and $u^\prime$ be its conjugate exponent. Then for any sequence $\{c_\alpha\}_{|\alpha|=k}$ of complex numbers, there exists a choice of signs $\{\epsilon_\alpha\}_{|\alpha|=k}$ such that
\begin{align*}
\sup_{z\in B(Y)}\left|\sum_{|\alpha|=k}\epsilon_\alpha c_\alpha z^\alpha\right|
\leq 
\begin{cases}
C(\log (n\log k))\sup\limits_{|\alpha|=k}\left\{|c_\alpha|\left(\frac{\alpha!}{k!}\right)\right\}\sup\limits_{\|z\|_Y\leq 1}\left(\sum_{m=1}^n|z_m|\right)^{k}\,\,\mbox{if}\,\, u=1,\\
C(n\log k)^{\frac{1}{u^\prime}}\sup\limits_{|\alpha|=k}\left\{|c_\alpha|\left(\frac{\alpha!}{k!}\right)^{\frac{1}{u}}\right\}\sup\limits_{\|z\|_Y\leq 1}\left(\sum_{m=1}^n|z_m|^u\right)^{\frac{k}{u}}\,\,\mbox{if}\,\, u\neq 1.
\end{cases}
\end{align*}

\vspace{3pt}
\noindent
Similarly, the following consequence of a famous result by Maurey and Pisier yields the desired upper bound when $\dim(X)=\infty$.

\vspace{3pt}
\noindent
\textbf{Theorem E (see f.i. \cite[Theorem 14.5]{Def3}).}
Given any complex Banach space $X$ with $\mbox{dim}(X)=\infty$, there exist $x_1, x_2,\cdots, x_n\in X$ for each $n\in\IN$ such that 
$$
\frac{1}{2}\|z\|_\infty\leq\left\|\sum_{k=1}^n x_kz_k\right\|\leq\|z\|_{\textrm{Cot}(X)}
$$
for every choice of $z=(z_1, z_2,\cdots, z_n)\in\IC^n$. Clearly, setting $z=e_k$-$e_k$ being the $k$-th canonical basis vector of $\IC^n$-gives $\|x_k\|\geq 1/2$.

\vspace{3pt}
\noindent
Furthermore, we state the following proposition, which can be established by adopting exactly similar lines of argument as in the proof of \cite[Proposition 19.2]{Def3}. We choose to include a proof for the sake of completeness.
\bprop\label{N3prop1}
For any complete Reinhardt domain $R\subset\IC^n$, we have
$$K_n^p(H^\infty(R, X))\geq K_n^p(H^\infty(\ID^n, X)).$$
\eprop
\bpf
For $\sigma=(r_1, r_2,\cdots, r_n)$ with $r_i>0$ for all $1\leq i\leq n$, we equip $\IC^n$ with the norm $\|z\|:=\sup_{1\le i\leq n}|z_i|/r_i$ and denote the resulting Banach space by $l_\infty^n(\sigma)$. The unit ball of this space is denoted by $B(l_\infty^n(\sigma))$. By the definition of complete Reinhardt domains, it is known that
$$
R=\bigcup_{B(l_\infty^n(\sigma))\subset R}B(l_\infty^n(\sigma)).
$$
Now, for $f(z)=\sum_{k=0}^\infty\sum_{|\alpha|=k}x_\alpha z^\alpha\in H^\infty(R, X)$ with $\|f\|_\infty\leq 1$, we form 
$$f_\sigma(z):=\sum_{k=0}^\infty\sum_{|\alpha|=k}x_\alpha \sigma^\alpha z^\alpha=\sum_{k=0}^\infty\sum_{|\alpha|=k}y_\alpha z^\alpha\in H^\infty(\ID^n, X)$$ for each $\sigma$ such that $B(l_\infty^n(\sigma))\subset R$. Evidently, $\|f_\sigma\|_\infty\leq 1$. Thus, for all $z\in K_n^p(H^\infty(\ID^n, X))\ID^n$, we have
$\sum_{k=1}^\infty\sum_{|\alpha|=k}\|y_\alpha z^\alpha\|^p\leq 1$.
But $y_\alpha=\sigma^\alpha x_\alpha$, and hence 
$\sum_{k=1}^\infty\sum_{|\alpha|=k}\|x_\alpha z^\alpha\|^p\leq 1$
for all $z\in K_n^p(H^\infty(\ID^n, X))B(l_\infty^n(\sigma))$. Due to the fact that $R$ is union of such \lq$B(l_\infty^n(\sigma))$' s, our proof is complete.
\epf

Let us conclude this section with a chain of inequalities, which will be used multiple times in the forthcoming proofs.
\be\label{N3eq20}
\sum_{|\alpha|=k}1={n+k-1\choose k}\leq e^k\left(1+\frac{n}{k}\right)^k\leq (2e)^k\max\left\{1, \left(\frac{n}{k}\right)^k\right\}.
\ee
\section{Proofs of the theorems}
\bpf[Proof of Theorem \ref{N3thm1}]
We start by observing the (possibly well-known) fact that for any $f(z)=\sum_{k=0}^\infty\sum_{|\alpha|=k}x_\alpha z^\alpha\in H^s(\ID^n, X)$ with $\|f\|_s\leq 1$: 
\be\label{N3eq1}
\left(\bigintsss_{\mathbb{T}^n}\left\|\sum_{|\alpha|=k}x_\alpha z^\alpha\right\|^s dz\right)^{\frac{1}{s}}\leq 1\,\mbox{for}\,s<\infty,\,\sup_{z\in\mathbb{T}^n}\left\|\sum_{|\alpha|=k}x_\alpha z^\alpha\right\|\leq 1\,\mbox{for}\,s=\infty.
\ee
Indeed, for any $\rho\in(0,1)$, we have for $s<\infty$:
\begin{align*}
\left(\bigintsss_{\mathbb{T}^n}\left\|\sum_{|\alpha|=k}x_\alpha z^\alpha\right\|^s dz\right)^{\frac{1}{s}}&=\left(\bigintsss_{\mathbb{T}^n}\left\|\frac{1}{2\pi}\bigintsss_{\theta=0}^{2\pi}\frac{f(\rho e^{i\theta}z)}{\rho^ke^{ik\theta}}d\theta\right\|^s dz\right)^{\frac{1}{s}}\\
&\leq\frac{1}{\rho^k}\left(\frac{1}{2\pi}\int_{\mathbb{T}^n}\int_{\theta=0}^{2\pi}\left\|f(\rho e^{i\theta}z)\right\|^s d\theta dz\right)^{\frac{1}{s}}\\
&=\frac{1}{\rho^k}\left(\frac{1}{2\pi}\int_{\theta=0}^{2\pi}\int_{\mathbb{T}^n}\left\|f(\rho z)\right\|^s dz d\theta \right)^{\frac{1}{s}}\leq\frac{\|f\|_s}{\rho^k},
\end{align*}
and thus by letting $\rho\to 1-$ we get the first inequality in $(\ref{N3eq1})$. For $s=\infty$, we consider the function $g(u)=f(uz)=x_0+\sum_{k=1}^\infty\left(\sum_{|\alpha|=k}x_\alpha z^\alpha\right)u^k$ for $u\in\ID$ and for any fixed $z\in\mathbb{T}^n$. $g$ is now a holomorphic function from $\ID$ to $X$ with $k$-th coefficient $\sum_{|\alpha|=k}x_\alpha z^\alpha$, and $\|g(u)\|\leq 1$ for all $u\in\ID$. As a result, for any $\phi\in X^*$ with $\|\phi\|\leq 1$, $\phi\circ g:\ID\to\ov{\ID}$ is holomorphic with $k$-th coefficient $\phi\left(\sum_{|\alpha|=k}x_\alpha z^\alpha\right)$. Therefore
$\left|\phi\left(\sum_{|\alpha|=k}x_\alpha z^\alpha\right)\right|\leq 1$, which in turn implies
$\sup_{z\in\mathbb{T}^n}\left\|\sum_{|\alpha|=k}x_\alpha z^\alpha\right\|\leq 1$.

With this information in hand, we are now ready to start proving the part (i) of this theorem.

\vspace{3pt}
\noindent
\underline{Proof of part (i)}: For any complex Banach space $X$ (finite/infinite dimensional) with a finite cotype $q$, it is known from \cite[Theorem 3.1]{Car} that there exists $C_1>0$ such that
$$
\left(\sum_{|\alpha|=k}\|x_\alpha\|^q\right)^{\frac{1}{q}}\leq C_1^k\left(\bigintsss_{\mathbb{T}^n}\left\|\sum_{|\alpha|=k}x_\alpha z^\alpha\right\|^q dz\right)^{\frac{1}{q}}.
$$
It is easy to see that 
$$
\left(\bigintsss_{\mathbb{T}^n}\left\|\sum_{|\alpha|=k}x_\alpha z^\alpha\right\|^q dz\right)^{\frac{1}{q}}\leq \sup_{z\in\mathbb{T}^n}\left\|\sum_{|\alpha|=k}x_\alpha z^\alpha\right\|,
$$
and for any $s<\infty$ the polynomial Kahane inequality (cf. \cite[Proposition 1.2]{Car1}) guarantees the existence of $C_2>0$ such that
$$
\left(\bigintsss_{\mathbb{T}^n}\left\|\sum_{|\alpha|=k}x_\alpha z^\alpha\right\|^q dz\right)^{\frac{1}{q}}\leq C_2^k \left(\bigintsss_{\mathbb{T}^n}\left\|\sum_{|\alpha|=k}x_\alpha z^\alpha\right\|^s dz\right)^{\frac{1}{s}}.
$$
The above-mentioned inequalities and $(\ref{N3eq1})$ ensure the existence of $C>0$ such that
\be\label{N3eq2}
\left(\sum_{|\alpha|=k}\|x_\alpha\|^q\right)^{\frac{1}{q}}\leq C^k.
\ee
Now, for any $p>\textrm{Cot}(X)$, there exists $q\in[\textrm{Cot}(X), p)$ such that $X$ has cotype $q$, and hence $X$ has cotype $p$. From $(\ref{N3eq2})$ we immediately get $\sum_{|\alpha|=k}\|x_\alpha\|^p\leq C^{pk}$.
Therefore $$\sum_{k=1}^\infty r^{pk}\sum_{|\alpha|=k}\|x_\alpha\|^p\leq \sum_{k=1}^\infty (Cr)^{pk}\leq 1$$
for all $r\leq 1/(2^{1/p}C)$, which is independent of $n$. On the other hand, it is evident that $K_n^p(H^s(\ID^n, X))\leq K_1^p(H^s(\ID, \IC))\leq 1$. As a result, $K_n^p(H^s(\ID^n, X))\sim 1$ for $p>\textrm{Cot}(X)$.

For $p\leq\textrm{Cot}(X)$, $p\leq q$ whenever $X$ has the cotype $q$. Using the inequalities $(\ref{N3eq20})$ and $(\ref{N3eq2})$ appropriately, we get 
\begin{align*}
\sum_{k=1}^\infty r^{pk}\sum_{|\alpha|=k}\|x_\alpha\|^p
&\leq \sum_{k=1}^\infty r^{pk}\left(\sum_{|\alpha|=k}\|x_\alpha\|^q\right)^{\frac{p}{q}}\left(\sum_{|\alpha|=k} 1\right)^{1-\frac{p}{q}}\\
&\leq \sum_{k=1}^\infty (e^{\frac{1}{p}-\frac{1}{q}}Cr)^{pk}(2n)^{\left(1-\frac{p}{q}\right)k}=\sum_{k=1}^\infty \left(C^\prime r n^{\frac{1}{p}-\frac{1}{q}}\right)^{pk},
\end{align*}
for some constant $C^\prime>0$ which is independent of $n$, and hence the above quantity is less than or equal to $1$ for all $r\leq 1/(2^{\frac{1}{p}}C^\prime n^{\frac{1}{p}-\frac{1}{q}})$. Therefore, in this case 
\be\label{N3eq3}
K_n^p(H^s(\ID^n, X))\geq C^{\prime\prime}\left(\frac{1}{n}\right)^{\frac{1}{p}-\frac{1}{q}}
\ee
for some constant $C^{\prime\prime}>0$. So far, all the results obtained by us remain valid for both finite and infinite dimensional Banach spaces $X$. Assuming now $\mbox{dim}(X)=\infty$, consider the $X$-valued holomorphic function $F(z)=\sum_{k=1}^n x_kz_k$ on $\ID^n$, where $z=(z_1, z_2, \cdots, z_n)\in\ID^n$ and \lq$x_k$' s are as in Theorem E. The definition of $K_n^p(H^s(\ID^n, X))$ combined with Theorem E gives 
$$
\frac{n}{2^p}\leq\sum_{k=1}^n\|x_k\|^p\leq\frac{\left\|\sum_{k=1}^n x_kz_k\right\|_s^p}{(K_n^p(H^s(\ID^n, X)))^p}\leq\frac{\sup_{z\in\ID^n}\|z\|_{\textrm{Cot}(X)}^p}{(K_n^p(H^s(\ID^n, X)))^p}\leq \frac{n^{\frac{p}{\textrm{Cot}(X)}}}{(K_n^p(H^s(\ID^n, X)))^p},
$$
i.e.
\be\label{N3eq4}
K_n^p(H^s(\ID^n, X))\leq 2\left(\frac{1}{n}\right)^{\frac{1}{p}-\frac{1}{\textrm{Cot}(X)}}.
\ee
The inequalities $(\ref{N3eq3})$ and $(\ref{N3eq4})$ together finish the proof for infinite dimensional Banach spaces $X$ with finite cotype. 

If $X$ does not have a finite cotype, i.e. $\textrm{Cot}(X)=\infty$, then again $(\ref{N3eq4})$ readily gives $K_n^p(H^s(\ID^n, X))\leq 2\left(1/n\right)^{1/p}$, and the required lower bound  $K_n^p(H^s(\ID^n, X))\succ (1/n)^{1/p}$ follows from the fact that $\|x_\alpha\|\leq 1$ (can be derived immediately from the Cauchy integral formula) and a computation exactly similar to the one used in proving $(\ref{N3eq3})$. Our proof for the first part of this theorem is therefore complete.

\vspace{3pt}
\noindent
\underline{Proof of part (ii)}: For finite dimensional $X$, it is known that $X$ has cotype $2$.  From part (i), we thus immediately conclude $K_n^p(H^s(\ID^n, X))\sim 1$ for all $s\geq 1$ and $p\geq\textrm{Cot}(X)=2$. We also get from $(\ref{N3eq3})$ that for $p<2$,
\be\label{N3eq5}
K_n^p(H^s(\ID^n, X))\succ\left(\frac{1}{n}\right)^{\frac{1}{p}-\frac{1}{2}}.
\ee
Set now 
\be\label{N3eq7}
G(z):=\sum_{k=1}^\infty\sum_{|\alpha|=k}r_0^{\frac{pk}{2-p}}z^\alpha,\,\,\mbox{where}\,\,r_0:=\left(1-\left(\frac{1}{2}\right)^{\frac{1}{n}}\right)^{\frac{1}{p}-\frac{1}{2}}.
\ee
Note that $\|G\|_s\leq\|G\|_2=1$ for $s\in[1, 2]$, and 
$\sum_{k=1}^\infty\sum_{|\alpha|=k}\left(r_0^{\frac{p}{2-p}}r\right)^{pk}\leq 1$
if and only if $r\leq r_0$. This means 
\be\label{N3eq6}
K_n^p(H^s(\ID^n, X))\leq K_n^p(H^s(\ID^n, \IC))\leq r_0.
\ee
A little calculation reveals that $\lim_{n\to\infty}n^{\frac{1}{p}-\frac{1}{2}}r_0=(\ln (2))^{\frac{1}{p}-\frac{1}{2}}$, which thereby gives $K_n^p(H^s(\ID^n, X))\prec (1/n)^{\frac{1}{p}-\frac{1}{2}}$. Together with $(\ref{N3eq5})$, this yields
$$
K_n^p(H^s(\ID^n, X))\sim \left(\frac{1}{n}\right)^{\frac{1}{p}-\frac{1}{2}}
$$
whenever $p<2$ and $1\leq s\leq 2$. 

We now take care of the $s>2$ part of this proof. For any $k$-homogeneous complex polynomial $P(z)=\sum_{|\alpha|=k}c_\alpha z^\alpha$, we have
$$
\left(\sum_{|\alpha|=k}|c_\alpha|^2\right)^{\frac{1}{2}}
=\left(\int_{\mathbb{T}^n} |P(z)|^2 dz\right)^{\frac{1}{2}},
$$
and due to the hypercontractive polynomial Bohnenblust-Hille inequality (see \cite[Theorem 1]{Seip}), there exists a constant $C_1>0$ such that 
$$
\left(\sum_{|\alpha|=k}|c_\alpha|^{\frac{2k}{k+1}}\right)^{\frac{k+1}{2k}}
\leq (C_1)^k\sup_{z\in\mathbb{T}^n}\left|P(z)\right|.
$$
In Theorem A, we now take $p_0=q_0=2$, $p_1=\infty, q_1=2k/(k+1)$, and $p_2=s$. As consequences of the above two inequalities, we have $M_0\leq 1, M_1\leq (C_1)^k$. Moreover
$$
\theta=1-\frac{2}{s},\,\,\frac{1}{q_2}=\frac{1}{2}+\frac{1}{2k}\left(1-\frac{2}{s}\right),
$$
and Theorem A further asserts that $M_2\leq(C_2)^k$ ($C_2=C_1^\theta$). This is same as saying
\be\label{N3eq8}
\left(\sum_{|\alpha|=k}|c_\alpha|^{q_2}\right)^{\frac{1}{q_2}}\leq (C_2)^k\times 
\begin{cases}
\left(\int_{\mathbb{T}^n} |P(z)|^s dz\right)^{\frac{1}{s}}\,\,\mbox{for}\,\, s<\infty,\\
\sup_{z\in\mathbb{T}^n}|P(z)|\,\,\mbox{for}\,\, s=\infty.
\end{cases}
\ee
Now for any $f(z)=\sum_{k=0}^\infty\sum_{|\alpha|=k}x_\alpha z^\alpha\in H^s(\ID^n, X)$ with $\|f\|_s\leq 1$, we construct
\be\label{N3eq12}
f_1(z):=\phi(f(z))=\sum_{k=0}^\infty\phi\left(\sum_{|\alpha|=k}x_\alpha z^\alpha\right)=\sum_{k=0}^\infty\sum_{|\alpha|=k}\phi(x_\alpha)z^\alpha,
\ee
where $\phi\in B(X^*)$-the open unit ball in $X^*$. It is easy to see that $f_1\in H^s(\ID^n, \IC)$ with $\|f_1\|_s\leq 1$ for all $s\in[1, \infty]$, and from $(\ref{N3eq1})$ it could also be readily verified that 
$$
\left(\bigintsss_{\mathbb{T}^n}\left|\phi\left(\sum_{|\alpha|=k}x_\alpha z^\alpha\right)\right|^s dz\right)^{\frac{1}{s}}\leq 1\,\mbox{for}\,s<\infty,\,\sup_{z\in\mathbb{T}^n}\left|\phi\left(\sum_{|\alpha|=k}x_\alpha z^\alpha\right)\right|\leq 1\,\mbox{for}\,s=\infty.
$$
As $X$ is finite dimensional, $\mbox{dim}(I(X))<\infty$ in this case, where $I:X\to X$ is the identity operator. Thus using Theorem B, we have $I\in\Pi_t(X, X)$ for all $t\geq 1$. Therefore
$$
\left(\sum_{|\alpha|=k}\|x_\alpha\|^{q_2}\right)^{\frac{1}{q_2}}\leq \pi_{q_2}(I)\sup_{\phi\in B(X^*)}\left(\sum_{|\alpha|=k}|\phi(x_\alpha)|^{q_2}\right)^{\frac{1}{q_2}}
$$
for all $k\in\IN$. For $s<\infty$, $q_2>1$ for all $k\geq 1$, and for $s=\infty$, $q_2=1$ for $k=1$ and $q_2>1$ otherwise. Hence the Theorem C asserts that $\pi_{q_2}(I)\leq\pi_1(I)$, i.e. there exists a constant $D=\pi_1(I)$ (depending only on $X$) such that
$$
\left(\sum_{|\alpha|=k}\|x_\alpha\|^{q_2}\right)^{\frac{1}{q_2}}\leq D\sup_{\phi\in B(X^*)}\left(\sum_{|\alpha|=k}|\phi(x_\alpha)|^{q_2}\right)^{\frac{1}{q_2}}.
$$
Combining $(\ref{N3eq8})$ with the above inequality, we obtain
$$
\left(\sum_{|\alpha|=k}\|x_\alpha\|^{q_2}\right)^{\frac{1}{q_2}}\leq D\sup_{\phi\in B(X^*)}(C_2)^k\times
\begin{cases}
\left(\bigintss_{\mathbb{T}^n}\left|\phi\left(\sum\limits_{|\alpha|=k}x_\alpha z^\alpha\right)\right|^s dz\right)^{\frac{1}{s}}\,\mbox{for}\,s<\infty,\\
\sup_{z\in\mathbb{T}^n}\left|\phi\left(\sum\limits_{|\alpha|=k}x_\alpha z^\alpha\right)\right|\,\mbox{for}\,s=\infty.
\end{cases}
$$
It is therefore evident that $$\left(\sum_{|\alpha|=k}\|x_\alpha\|^{q_2}\right)^{\frac{1}{q_2}}\leq D(C_2)^k.$$
For any given $p\in(0, 2)$ now, let us set
$k_0:=\left\lfloor p\left(1-(2/s)\right)/\left(2-p\right)\right\rfloor+1$ ($\lfloor.\rfloor$ is the floor function).
Clearly, for $k<k_0$ we have $q_2\leq p$, which in turn gives 
$$\left(\sum_{|\alpha|=k}\|x_\alpha\|^{p}\right)^{\frac{1}{p}}\leq \left(D^{\frac{1}{k}}C_2\right)^k\leq(C_2^\prime)^k,$$
$C_2^\prime>0$ being a new constant. For $k\geq k_0$, we have $q_2\geq p$, and hence
$$
\sum_{|\alpha|=k}\|x_\alpha\|^{p}\leq\left(\sum_{|\alpha|=k}\|x_\alpha\|^{q_2}\right)^{\frac{p}{q_2}}\left(\sum_{|\alpha|=k} 1\right)^{1-\frac{p}{q_2}}
\leq (C_3)^{k}\left(\max\left\{1, \left(\frac{n}{k}\right)^{1-\frac{p}{q_2}}\right\}\right)^k
$$
(see (\ref{N3eq20})) for some new constant $C_3>0$.
A little computation reveals that the function 
$g(x)=x^{1-\frac{p}{2}}n^{\frac{p}{2x}(1-\frac{2}{s})}:(0, \infty)\to (0, \infty)$ satisfies
$$
g(x)\geq g\left(\frac{p\left(1-\frac{2}{s}\right)}{2-p}\log n\right)=C_4(\log n)^{1-\frac{p}{2}}
$$
for some constant $C_4>0$.
Using this and the fact that $k^{1/k}\leq e^{1/e}$ for all $k\geq 1$, we observe that
$$
\left(\frac{n}{k}\right)^{1-\frac{p}{q_2}}
=\frac{n^{1-\frac{p}{2}}k^{\frac{p}{2k}\left(1-\frac{2}{s}\right)}}{g(k)}
\leq C_5\left(\frac{n}{\log n}\right)^{1-\frac{p}{2}},
$$
$C_5>0$ being another constant, and therefore
\be\label{N3eq11}
\sum_{|\alpha|=k}\|x_\alpha\|^{p}\leq \left(C_6\left(\frac{n}{\log n}\right)^{1-\frac{p}{2}}\right)^k
\ee
for a constant $C_6>0$. It is now clear that
\begin{align*}
\sum_{k=1}^\infty r^{pk}\sum_{|\alpha|=k}\|x_\alpha\|^p
&=\sum_{k=1}^{k_0-1} r^{pk}\sum_{|\alpha|=k}\|x_\alpha\|^p+\sum_{k=k_0}^\infty r^{pk}\sum_{|\alpha|=k}\|x_\alpha\|^p\\
&\leq\sum_{k=1}^{k_0-1}r^{pk}(C_2^\prime)^{pk}+\sum_{k=k_0}^\infty r^{pk}\left(C_6\left(\frac{n}{\log n}\right)^{1-\frac{p}{2}}\right)^k,
\end{align*}
which is less than or equal to $1$ if $r=C((\log n)/n)^{(1/p)-(1/2)}$, $C>0$ is a sufficiently small constant. That is to say 
$$
K_n^p(H^s(\ID^n, X))\succ\left(\frac{\log n}{n}\right)^{\frac{1}{p}-\frac{1}{2}}.
$$
Also, using exactly the same argument as in \cite[p. 76]{Dja}, it is possible to establish that 
$$
K_n^p(H^\infty(\ID^n, \IC))\prec\left(\frac{\log n}{n}\right)^{\frac{1}{p}-\frac{1}{2}}
$$
for all $p\in(0, 2)$. Since $K_n^p(H^s(\ID^n, X))\leq K_n^p(H^s(\ID^n, \IC))\leq K_n^p(H^\infty(\ID^n, \IC))$,
this step finishes our proof for the $s>2$ case.

\vspace{3pt}
\noindent
\underline{Proof of part (iii)}: For any $f(z)=\sum_{k=0}^\infty\sum_{|\alpha|=k}x_\alpha z^\alpha\in H^2(\ID^n, \IC)$ with $\|f\|_2\leq 1$, it is easy to see that $\sum_{k=0}^\infty\sum_{|\alpha|=k}|x_\alpha|^2\leq 1 \,(\lq x_\alpha \textrm{' s are complex numbers now})$. Therefore, for $p<2$:
\begin{align*}
\sum_{k=1}^\infty r^{pk}\sum_{|\alpha|=k}|x_\alpha|^p&\leq\left(\sum_{k=1}^\infty\sum_{|\alpha|=k}|x_\alpha|^2\right)^{\frac{p}{2}}\left(\sum_{k=1}^\infty\sum_{|\alpha|=k}r^{\frac{2pk}{2-p}}\right)^{1-\frac{p}{2}}\\
&\leq \left(\sum_{k=1}^\infty{n+k-1 \choose k}r^{\frac{2pk}{2-p}}\right)^{1-\frac{p}{2}}
=\left(\frac{1}{\left(1-r^{\frac{2p}{2-p}}\right)^n}-1\right)^{1-\frac{p}{2}},
\end{align*}
which is less than or equal to $1$ for $r\leq r_0$, $r_0$ as defined in $(\ref{N3eq7})$. Hence we conclude that $K_n^p(H^2(\ID^n, \IC))\geq r_0$, and from $(\ref{N3eq6})$ we already have $K_n^p(H^2(\ID^n, \IC))\leq r_0$. This completes the proof for part (iii) of this theorem.
\epf
\bpf[Proof of Theorem \ref{N3thm2}]
Since the case $t=\infty$ is already covered in Theorem \ref{N3thm1}, we only need to consider $t\in[1, \infty)$.
Let us begin with a coefficient estimate for $f(z)=\sum_{k=0}^\infty\sum_{|\alpha|=k}x_\alpha z^\alpha\in H^\infty(B(l_t^n), X)$ with $\|f(z)\|\leq 1$ for all $z\in B(l_t^n)$. We construct $f_1$ from $f$ as in $(\ref{N3eq12})$. It is clear that $|f_1(z)|\leq 1$ for all $z\in B(l_t^n)$, and hence from the inequality $(3.7)$ of \cite[p. 24]{Def4} we conclude
$$
\|x_\alpha\|=\sup_{\phi\in B(X^*)}\left|\phi(x_\alpha)\right|\leq e^{\frac{k}{t}}\left(\frac{k!}{\alpha!}\right)^{\frac{1}{t}}
$$
for all multi-indices $\alpha$ with $|\alpha|=k$. This inequality plays a crucial role in the proof of the first part of this theorem.

\vspace{3pt}
\noindent
\underline{Proof of part (i)}: For any $z=(z_1, z_2, \cdots, z_n)\in B(l_t^n)$ and $r\in (0, 1)$, we have
\begin{align*}
\sum_{k=1}^\infty\sum_{|\alpha|=k}\|x_\alpha (rz)^\alpha\|^t
&\leq\sum_{k=1}^\infty r^{kt}e^k\sum_{|\alpha|=k}\left(\frac{k!}{\alpha!}\right)|z_1|^{\alpha_1t}|z_2|^{\alpha_2t}\cdots|z_n|^{\alpha_nt}\\
&=\sum_{k=1}^\infty r^{kt}e^k\left(|z_1|^t+|z_2|^t+\cdots+|z_n|^t\right)^k
\leq \sum_{k=1}^\infty r^{kt}e^k,
\end{align*}
which is less than or equal to $1$ if $r\leq (1/(2e))^{1/t}$, i.e. $K_n^p(H^\infty(B(l_t^n), X))\sim 1$ for $p\geq t$. Also, since $B(l_t^n)$ is a complete Reinhardt domain in $\IC^n$, using the part (i) of Theorem \ref{N3thm1} and Proposition \ref{N3prop1} gives $K_n^p(H^\infty(B(l_t^n), X))\sim 1$ for $p>\textrm{Cot}(X)$. Therefore, we get 
$$
K_n^p(H^\infty(B(l_t^n), X))\sim 1
$$
for $p>\min\{\textrm{Cot}(X), t\}$. When $p\leq\min\{\textrm{Cot}(X), t\}$, $p\leq q$ provided that $X$ has finite cotype $q$. Denote by $\gamma:=\min\{q, t\}$. Then it is clear from the above discussion that there exists $r_2>0$ independent of $n$, such that
\be\label{N2eq16}
\sum_{k=1}^\infty\sum_{|\alpha|=k}\|x_\alpha (r_2z)^\alpha\|^\gamma\leq 1,\,\,z\in B(l_t^n).
\ee
Therefore, for any $z=(z_1, z_2, \cdots, z_n)\in B(l_t^n)$ and $r=r_1r_2\in (0, 1)$, using the inequality $(\ref{N3eq20})$ gives
\begin{align*}
\sum_{k=1}^\infty\sum_{|\alpha|=k}\|x_\alpha (rz)^\alpha\|^p
&\leq\sum_{k=1}^\infty r_1^{pk}\left(\sum_{|\alpha|=k}\|x_\alpha (r_2z)^\alpha\|^\gamma\right)^{\frac{p}{\gamma}}\left(\sum_{|\alpha|=k} 1\right)^{1-\frac{p}{\gamma}}\\
&\leq \sum_{k=1}^\infty (e^{\frac{1}{p}-\frac{1}{\gamma}}r_1)^{pk}(2n)^{\left(1-\frac{p}{\gamma}\right)k}=\sum_{k=1}^\infty \left(C^\prime r_1n^{\frac{1}{p}-\frac{1}{\gamma}}\right)^{pk}
\end{align*}
($C^\prime>0$ being a constant), which is less than or equal to $1$ if $r_1=C(1/n)^{\frac{1}{p}-\frac{1}{\gamma}}$ for some sufficiently small $C>0$. This yields
\be\label{N3eq13}
K_n^p(H^\infty(B(l_t^n), X))\succ\left(\frac{1}{n}\right)^{\frac{1}{p}-\frac{1}{\gamma}}.
\ee
Let us note that all the results obtained upto this point remain valid for both finite and infinite dimensional Banach spaces. Now to obtain the upper estimate on $K_n^p(H^\infty(B(l_t^n), X))$, we need to assume $\dim (X)=\infty$. This portion of our proof proceeds exactly like the proof of $(\ref{N3eq4})$, i.e. we define $F(z)=\sum_{k=1}^n x_k z_k$ on $B(l_t^n)$, where \lq$x_k$' s are as in Theorem E. The definition of $K_n^p(H^\infty(B(l_t^n), X))$ allows us to set $$z=K_n^p(H^\infty(B(l_t^n), X))\left(\left(\frac{1}{n}\right)^{\frac{1}{t}}, \left(\frac{1}{n}\right)^{\frac{1}{t}}, \cdots, \left(\frac{1}{n}\right)^{\frac{1}{t}}\right),$$ 
which, after an application of Theorem E, generates the following chain of inequalities:
$$
\frac{n}{2^p}\leq\sum_{k=1}^n\|x_k\|^p
\leq n^{\frac{p}{t}}\frac{\sup_{z\in B(l_t^n)}\left\|\sum_{k=1}^n x_kz_k\right\|^p}{(K_n^p(H^\infty(B(l_t^n), X)))^p}
\leq n^{\frac{p}{t}}\frac{\sup_{z\in B(l_t^n)}\|z\|_{\textrm{Cot}(X)}^p}{(K_n^p(H^\infty(B(l_t^n), X)))^p}.
$$
For $t\leq\textrm{Cot}(X)$, $\sup_{z\in B(l_t^n)}\|z\|_{\textrm{Cot}(X)}\leq 1$, which gives
$$
K_n^p(H^\infty(B(l_t^n), X))\leq 2\left(\frac{1}{n}\right)^{\frac{1}{p}-\frac{1}{t}},
$$
and for $t>\textrm{Cot}(X)$, $\sup_{z\in B(l_t^n)}\|z\|_{\textrm{Cot}(X)}\leq n^{(1/\textrm{Cot}(X))-(1/t)}$, which gives
$$
K_n^p(H^\infty(B(l_t^n), X))\leq 2\left(\frac{1}{n}\right)^{\frac{1}{p}-\frac{1}{\textrm{Cot}(X)}}.
$$
This is same as saying 
\be\label{N3eq15}
K_n^p(H^\infty(B(l_t^n), X))\prec\left(\frac{1}{n}\right)^{\frac{1}{p}-\frac{1}{\min\{\textrm{Cot}(X), t\}}}.
\ee
The two inequalities $(\ref{N3eq13})$ and $(\ref{N3eq15})$ together complete the proof for infinite dimensional complex Banach spaces $X$ with finite cotype.

If $X$ possesses no finite cotype, the upper bound for $K_n^p(H^\infty(B(l_t^n), X))$ can be obtained directly from $(\ref{N3eq15})$ by putting $\textrm{Cot}(X)=\infty$. For the lower bound, we only have to repeat the proof for $(\ref{N3eq13})$ line-by-line, by replacing $\gamma$ with $t$. This finishes the first part of the proof of this theorem.

\vspace{3pt}
\noindent
\underline{Proof of part (ii)}: For $\dim(X)<\infty$, the assertion $K_n^p(H^\infty(B(l_t^n), X))\sim 1$ for $p\geq\min\{t, 2\}$ follows readily from part (i), as $X$  has cotype $2$ and $\textrm{Cot}(X)=2$. When $p<\min\{t, 2\}$, we again need to work with $f_1$ as in $(\ref{N3eq12})$, defined on $B(l_t^n)$. Now for any fixed $z\in\IC^n$ with $\|z\|_{l_t^n}\leq 1$, $f_1(uz):\ID\to\ov{\ID}$ is a holomorphic function of $u\in\ID$, with $k$-th coefficient $\phi\left(\sum_{|\alpha|=k}x_\alpha z^\alpha\right)$. As a consequence 
$$
\sup_{\|z\|_{l_t^n}\leq 1}\left|\phi\left(\sum_{|\alpha|=k}x_\alpha z^\alpha\right)\right|\leq 1.
$$
From \cite[Lemma 2.1]{Def5} and \cite[Theorem 1.2]{Def6}, we have for a complex $k$-homogeneous polynomial $P(z)=\sum_{|\alpha|=k}c_\alpha z^\alpha$ with $|P(z)|\leq 1, z\in B(l_t^n)$:
\be\label{N3eq18}
\sup_{\|z\|_{l_t^n}\leq 1}\sum_{|\alpha|=k}\left|c_\alpha z^\alpha\right|\leq\chi_M(\mathcal{P}({}^{k}l_t^n))
\leq C^k\left(1+\frac{n}{k}\right)^{(k-1)\left(1-\frac{1}{\delta}\right)},
\ee
where $\delta:=\min\{2, t\}$. We clarify that $\chi_M(\mathcal{P}({}^{k}X))$ is the unconditional basis constant associated with the basis consisting of the monomials $z^\alpha$, for the space $\mathcal{P}({}^{k}X)$ of continuous $k$-homogeneous scalar-valued polynomials $P$ on a given Banach space $X$. This space is equipped with the norm
$\|P\|=\sup_{\|x\|\leq 1}|P(x)|$ (see \cite{Def5} for more details).
Making use of Theorem B, we now get
\begin{align*}
\sum_{|\alpha|=k}\|x_\alpha z^\alpha\|\leq \pi_{1}(I)\sup_{\phi\in B(X^*)}\left(\sum_{|\alpha|=k}|\phi(x_\alpha)z^\alpha|\right)
&\leq \pi_{1}(I)\chi_M(\mathcal{P}({}^{k}l_t^n))\\
&\leq\pi_{1}(I)C^k\left(1+\frac{n}{k}\right)^{(k-1)\left(1-\frac{1}{\delta}\right)}
\end{align*}
for all $z\in B(l_t^n)$. It is therefore clear that by exactly the same calculations as in the proof of \cite[Theorem 1.1]{Def6} (see \cite[p. 145]{Def6} in particular), we obtain
$$
K_X(1,t):=K_n^1(H^\infty(B(l_t^n), X))\succ\left(\frac{\log n}{n}\right)^{\left(1-\frac{1}{\delta}\right)}.
$$
Also, inequality $(\ref{N2eq16})$ is now true for $\gamma=\delta$.
For any $p\in [1, \delta)$, we observe that
$(1/p)=(1-\theta)/1+(\theta/\delta)$
if $\theta=(1-(1/p))/(1-(1/\delta))$. Setting now $r=(K_X(1,t))^{1-\theta}r_2^\theta$ ($r_2$ being the $n$-independent constant for which $(\ref{N2eq16})$ is satisfied, with $\gamma$ replaced by $\delta$), we get
\begin{align*}
\sum_{k=1}^\infty\sum_{|\alpha|=k}\|x_\alpha (rz)^\alpha\|^p
&=\sum_{k=1}^\infty\sum_{|\alpha|=k}\|x_\alpha (K_X(1,t)z)^\alpha\|^{p(1-\theta)}\|x_\alpha (r_2z)^\alpha\|^{p\theta}\\
&\leq\left(\sum_{k=1}^\infty\sum_{|\alpha|=k}\|x_\alpha (K_X(1,t)z)^\alpha\|\right)^{p(1-\theta)}\hspace{-6pt}\left(\sum_{k=1}^\infty\sum_{|\alpha|=k}\|x_\alpha (r_2z)^\alpha\|^\delta\right)^{\frac{p\theta}{\delta}}\hspace{-6pt}\leq 1,
\end{align*}
which implies $K_n^p(H^\infty(B(l_t^n), X))\succ(K_X(1,t))^{1-\theta}$. This gives the required lower bound for $1\leq p<\min\{t, 2\}$. For $p\in(0,1)$:
$$
\sum_{|\alpha|=k}\|x_\alpha z^\alpha\|^p \leq\left(\sum_{|\alpha|=k}\|x_\alpha z^\alpha\|\right)^p\left(\sum_{|\alpha|=k} 1\right)^{1-p}
\leq C_1^k\left(\max\left\{1, \left(\frac{n}{k}\right)^{1-\frac{p}{\delta}-\frac{p}{k}\left(1-\frac{1}{\delta}\right)}\right\}\right)^k
$$
(see (\ref{N3eq20}) and the aforesaid estimate for $\sum_{|\alpha|=k}\|x_\alpha z^\alpha\|$) for a constant $C_1>0$. We observe that the function
$g(x)=x^{1-\frac{p}{\delta}}n^{\frac{p}{x}\left(1-\frac{1}{\delta}\right)}:(0, \infty)\to(0, \infty)$
satisfies
$$
g(x)\geq g\left(\frac{p\left(1-\frac{1}{\delta}\right)}{1-\frac{p}{\delta}}\log n\right)=C_2\left(\log n\right)^{1-\frac{p}{\delta}},
$$
and therefore
$$
\left(\frac{n}{k}\right)^{1-\frac{p}{\delta}-\frac{p}{k}\left(1-\frac{1}{\delta}\right)}
=\frac{n^{1-\frac{p}{\delta}}k^{\frac{p}{k}\left(1-\frac{1}{\delta}\right)}}{g(k)}
\leq C_3\left(\frac{n}{\log n}\right)^{1-\frac{p}{\delta}},
$$
where $C_2$ and $C_3$ are two new constants.
It is now evident that
$$
\sum_{k=1}^\infty\sum_{|\alpha|=k}\|x_\alpha (rz)^\alpha\|^p
\leq \sum_{k=1}^\infty r^{kp}\left(C_3^\prime\left(\frac{n}{\log n}\right)^{1-\frac{p}{\delta}}\right)^k
$$
($C_3^\prime$ being another constant), which is less than or equal to $1$ if we choose $r=C((\log n)/n)^{(1/p)-(1/\delta)}$ for some sufficiently small $C>0$. We have thus established that
\be\label{N3eq17}
K_n^p(H^\infty(B(l_t^n), X))\succ\left(\frac{\log n}{n}\right)^{\frac{1}{p}-\frac{1}{\delta}}
\ee
for all $p$ satisfying $0<p<\delta=\min\left\{t, 2\right\}$.

It remains to obtain the upper bounds. Clearly, it is enough to derive these bounds for $X=\IC$. The argument will be divided in two parts: for the space $l_1^n$ and for all other spaces $l_t^n$ s. An essential tool of us is Stirling's formula 
$$
\lim_{k\to\infty}k!\left(\sqrt{2\pi k}(k/e)^k\right)^{-1}=1,
$$
which we will use frequently. For brevity, let us denote $K_n^p(H^\infty(B(l_t^n), \IC))$ by only $K(p,t)$ for the rest of this paper. For $t=1$, we already saw that $K(p,1)\sim 1$ for $p\geq 1$, so it is sufficient to consider $p<1$ in this situation. Set $c_\alpha=(k!/\alpha!)^{1/p}$ in the $u=1$ part of Theorem D. Using the definition of $p$-Bohr radius gives
\begin{align*}
&\sum_{|\alpha|=k}\left(\frac{k!}{\alpha!}\right)\left(\frac{K(p, 1)}{n}\right)^{pk}
\leq C^p\left(\log(n\log k)\right)^p\sup_{|\alpha|=k}\left\{\left(\frac{k!}{\alpha!}\right)^{1-p}\right\}\sup_{\|z\|_{l_1^n}\leq 1}\left(\sum_{m=1}^n|z_m|\right)^{pk}\\
&\implies\left(n^{1-p}(K(p,1))^p\right)^k\leq C_1\left(\log(n\log k)\right)^p\left(\frac{k^{k+\frac{1}{2}}}{e^k}\right)^{1-p}\\
&\implies\left(\frac{n}{k}\right)^{1-p}(K(p, 1))^p
\leq C_1^{\frac{1}{k}}\left(\log(n\log k)\right)^{\frac{p}{k}}\left(\frac{k^{\frac{1}{2k}}}{e}\right)^{1-p},
\end{align*}
where $C, C_1>0$ are constants. We now set $k=\lfloor \log n\rfloor$ and observe that
\begin{align*}
\limsup_{n\to\infty}\left(\frac{n}{\log n}\right)^{1-p}(K(p, 1))^p
&\leq\limsup_{n\to\infty}\left(\frac{n}{\lfloor \log n\rfloor}\right)^{1-p}(K(p, 1))^p\\
&\leq e^{p-1}\lim_{n\to\infty}\left(C_1^{\frac{1}{p}}\log(n\log\lfloor \log n\rfloor)\right)^{\frac{p}{\lfloor \log n\rfloor}}\lfloor \log n\rfloor^{\frac{1-p}{2\lfloor \log n\rfloor}}\\
&=e^{p-1},
\end{align*}
which implies $K(p, 1)\prec\left((\log n)/n\right)^{(1/p)-1}$. The proof for $l_1^n$ is therefore complete.

For $l_t^n$ with $t\neq 1$, the argument is similar but we need to do the proof separately for $t\in(1,2)$ and $t\in [2, \infty)$. When $t\in(1, 2)$, we set $u=t$ in the $u\neq 1$ part of Theorem D, and choose $c_\alpha=(k!/\alpha!)^{1/p}$ again. The definition of the $p$-Bohr radius gives
\begin{align*}
&\sum_{|\alpha|=k}\left(\frac{k!}{\alpha!}\right)\left(\frac{K(p, t)}{n^{\frac{1}{t}}}\right)^{pk}
\leq C^p\left(\log k\right)^{\frac{p}{t^\prime}}n^{\frac{p}{t^\prime}}\sup_{|\alpha|=k}\left\{\left(\frac{k!}{\alpha!}\right)^{1-\frac{p}{t}}\right\}\sup_{\|z\|_{l_t^n}\leq 1}\left(\sum_{m=1}^n|z_m|^t\right)^{\frac{pk}{t}}\\
&\implies\left(n^{1-\frac{p}{t}}(K(p, t))^p\right)^k\leq C_1\left(\log k\right)^{\frac{p}{t^\prime}}n^{\frac{p}{t^\prime}}\left(\frac{k^{k+\frac{1}{2}}}{e^k}\right)^{1-\frac{p}{t}}\\
&\implies\left(\frac{n}{k}\right)^{1-\frac{p}{t}}(K(p, t))^p
\leq C_1^{\frac{1}{k}}\left(\log k\right)^{\frac{p}{kt^\prime}}n^{\frac{p}{kt^\prime}}\left(\frac{k^{\frac{1}{2k}}}{e}\right)^{1-\frac{p}{t}},
\end{align*}
where $C, C_1>0$ are constants (not necessarily same with the previously used ones). Again, putting $k=\lfloor \log n\rfloor$ and by a little computation, we obtain
\begin{align*}
\limsup_{n\to\infty}\left(\frac{n}{\log n}\right)^{1-\frac{p}{t}}(K(p, t))^p
&\leq e^{\frac{p-t}{t}}\lim_{n\to\infty}\left(C_1^{\frac{t^\prime}{p}}\log \lfloor \log n\rfloor\right)^{\frac{p}{t^\prime\lfloor \log n\rfloor}}n^{\frac{p}{t^\prime\lfloor \log n\rfloor}}\lfloor \log n\rfloor^{\frac{1-\frac{p}{t}}{2\lfloor \log n\rfloor}}\\
&=e^{p-1},
\end{align*}
which gives us $K(p, t)\prec((\log n)/n)^{(1/p)-(1/t)},\, t\in(1,2)$. For the case $t\geq 2$, we have to choose $u=2$ in Theorem D, and by an argument exactly same as above we will get
\begin{align*}
&\left(n^{1-\frac{p}{t}}(K(p, t))^p\right)^k\leq C_1\left(\log k\right)^{\frac{p}{2}}n^{\frac{p}{2}}\left(\frac{k^{k+\frac{1}{2}}}{e^k}\right)^{1-\frac{p}{2}}\sup_{\|z\|_{l_t^n}\leq 1}\left(\sum_{m=1}^n|z_m|^2\right)^{\frac{pk}{2}}\\
&\implies\left(\frac{n}{k}\right)^{1-\frac{p}{2}}(K(p, t))^p
\leq C_1^{\frac{1}{k}}\left(\log k\right)^{\frac{p}{2k}}n^{\frac{p}{2k}}\left(\frac{k^{\frac{1}{2k}}}{e}\right)^{1-\frac{p}{2}}.
\end{align*}
From this step, it is clear that the method for the proof of the previous case (i.e. $t\in(1, 2)$) can be further adopted to arrive at the conclusion $K(p, t)\prec((\log n)/n)^{(1/p)-(1/2)}$. This completes our proof.
\epf
\bpf[Proof of Theorem \ref{N3thm3}]For a pluriharmonic function $f(z)=\sum_{k=0}^\infty\sum_{|\alpha|=k}a_\alpha z^\alpha+\ov{\sum_{k=1}^\infty\sum_{|\alpha|=k}b_\alpha z^\alpha}$ from $B(l_t^n)$ to $\ov{\ID}$, it is clear that both of the functions
$\mbox{Re}(f(z))$ and $\mbox{Im}(f(z))$
map $B(l_t^n)$ inside $\ov{\ID}$. As $\mbox{Im}(f(z))=\mbox{Re}(-if(z))$, and we can assume $f(0)\geq 0$ without loss of generality, it is easy to observe from \cite[Theorem 2.7]{Chen} that
$$
\sup_{\|z\|_{l_t^n}\leq 1}\left|\sum_{|\alpha|=k}\left(\frac{a_\alpha+b_\alpha}{2}\right)z^\alpha\right|\leq 1,\,\,\sup_{\|z\|_{l_t^n}\leq 1}\left|\sum_{|\alpha|=k}\left(\frac{a_\alpha-b_\alpha}{2}\right)z^\alpha\right|\leq 1.
$$
Our definition of $p$-Bohr radius for holomorphic functions readily yields
$$
\sum_{|\alpha|=k}\left|\left(\frac{a_\alpha+b_\alpha}{2}\right)z^\alpha\right|^p\leq\frac{1}{K(p,t)^{pk}}\,\,\textrm{and}\,\,
\sum_{|\alpha|=k}\left|\left(\frac{a_\alpha-b_\alpha}{2}\right)z^\alpha\right|^p\leq\frac{1}{K(p,t)^{pk}}
$$
for all $z\in B(l_t^n)$, $K(p,t)=K_n^p(H^\infty(B(l_t^n), \IC))$ as defined in the final stage of the proof of part (ii) of Theorem \ref{N3thm2}. Now for $p\geq 1$, a little computation reveals that
$$
\left(\sum_{|\alpha|=k}|a_\alpha z^\alpha|^p\right)^{\frac{1}{p}}\hspace{-5pt}\leq
\left(\sum_{|\alpha|=k}\left|\left(\frac{a_\alpha+b_\alpha}{2}\right)z^\alpha\right|^p\right)^{\frac{1}{p}}\hspace{-5pt}+\left(\sum_{|\alpha|=k}\left|\left(\frac{a_\alpha-b_\alpha}{2}\right)z^\alpha\right|^p\right)^{\frac{1}{p}}
\hspace{-5pt}\leq\frac{2}{K(p,t)^{k}}
$$
for all $z\in B(l_t^n)$. Similarly
$\left(\sum_{|\alpha|=k}|b_\alpha z^\alpha|^p\right)^{1/p}\leq 2/(K(p,t))^{k}$
for all $z\in B(l_t^n)$. As a result, for $r\in(0,1)$
$$
\sum_{k=1}^\infty\sum_{|\alpha|=k}\left(|a_\alpha|^p+|b_\alpha|^p\right)|(rz)^{\alpha}|^p\leq 2^{p+1}\sum_{k=1}^\infty\left(\frac{r}{K(p,t)}\right)^{pk},
$$
which is less than or equal to $1$ if we choose $r=CK(p, t)$ for $C$ small enough. Moreover, it is clear that $R_p(B(l_t^n))\leq K(p,t)$ for all $p>0$, and hence we have established the following:
$$
R_p(B(l_t^n))\sim
\begin{cases}
\left(\frac{\log n}{n}\right)^{\frac{1}{p}-\frac{1}{\min\{t,2\}}}\,\,\emph{for}\,\, 1\leq p<\min\{t,2\},\,\emph{and}\\
1\,\,\emph{for}\,\,p\geq\min\{t,2\}.   
\end{cases}
$$
It remains to prove that $R_p(B(l_t^n))\succ
\left((\log n)/n\right)^{(1/p)-(1/\delta)}$ for $p\in(0, 1)$, $\delta=\min\{t,2\}$ as defined in the proof of part (ii) of Theorem \ref{N3thm2}. Indeed, inequality $(\ref{N3eq18})$ asserts that
$$
\sum_{|\alpha|=k}\left|\left(\frac{a_\alpha+b_\alpha}{2}\right)z^\alpha\right|\leq\chi_M(\mathcal{P}({}^{k}l_t^n))\,\,\mbox{and}\,\,
\sum_{|\alpha|=k}\left|\left(\frac{a_\alpha-b_\alpha}{2}\right)z^\alpha\right|\leq\chi_M(\mathcal{P}({}^{k}l_t^n)),
$$
and as a consequence
$$
\sum_{|\alpha|=k}|a_\alpha z^\alpha|\leq 2\chi_M(\mathcal{P}({}^{k}l_t^n))
\leq 2C^k\left(1+\frac{n}{k}\right)^{(k-1)\left(1-\frac{1}{\delta}\right)}
$$
for all $z\in B(l_t^n)$. The above inequality remains valid if the left hand side is replaced by $\sum_{|\alpha|=k}|b_\alpha z^\alpha|$ as well. Therefore, making use of $(\ref{N3eq20})$, we get
\begin{align*}
\sum_{|\alpha|=k}\left(|a_\alpha z^\alpha|^p+|b_\alpha z^\alpha|^p\right)
&\leq\left(\left(\sum_{|\alpha|=k}|a_\alpha z^\alpha|\right)^p+\left(\sum_{|\alpha|=k}|b_\alpha z^\alpha|\right)^p\right)\left(\sum_{|\alpha|=k} 1\right)^{1-p}\\
&\leq \left(C_1^k+C_2^k\right)\left(\max\left\{1, \left(\frac{n}{k}\right)^{1-\frac{p}{\delta}-\frac{p}{k}\left(1-\frac{1}{\delta}\right)}\right\}\right)^k\\
&\leq C^k\left(\max\left\{1, \left(\frac{n}{k}\right)^{1-\frac{p}{\delta}-\frac{p}{k}\left(1-\frac{1}{\delta}\right)}\right\}\right)^k
\end{align*}
for some constants $C, C_1, C_2>0$ and for all $z\in B(l_t^n)$. The proof will now proceed exactly like the proof of $(\ref{N3eq17})$ (the $p\in(0,1)$ part) and hence we omit the details.
\epf

\section*{Acknowledgements}
The author is thankful to Prof. Bappaditya Bhowmik for his kind help in obtaining a softcopy of \cite{Dja}. He is also indebted to Dr. Samya Kumar Ray for some enriching conversations.

\end{document}